\documentclass[11pt, a4paper]{amsart}
\usepackage{amssymb,latexsym,amsmath, graphicx} 
\usepackage{enumerate,enumitem,color,hyperref} 
\usepackage{mathrsfs}
\usepackage{calligra}
\usepackage[T1]{fontenc}

\input{xypic}
\xyoption{all}

\title[\bf Congruences]{Ramanujan's congruence primes}

\author{ \bf Ellise Parnoff \ \& \ A. Raghuram} 

\date{\today}      

\address{Department of Mathematics, 
Fordham University at Lincoln Center, 113 W 60th St, New~York, NY 10023, USA.} 

\email{eparnoff@fordham.edu, \ araghuram@fordham.edu}

\numberwithin{equation}{section}   

\newtheorem{thm}[equation]{Theorem}

\newtheorem{rem}[equation]{Remark}

\usepackage{bm}
\usepackage{upgreek}
\makeatletter
\newcommand{\bfgreek}[1]{\bm{\@nameuse{up#1}}}

\let\oldtocsection=\tocsection
\let\oldtocsubsection=\tocsubsection
\let\oldtocsubsubsection=\tocsubsubsection
\renewcommand{\tocsection}[2]{\hspace{0em}\oldtocsection{#1}{#2}}
\renewcommand{\tocsubsection}[2]{\hspace{1em}\oldtocsubsection{#1}{#2}}
\renewcommand{\tocsubsubsection}[2]{\hspace{2em}\oldtocsubsubsection{#1}{#2}}

\begin{document}

\begin{abstract}
Ramanujan showed that $\tau(p) \equiv p^{11}+1 \pmod{691}$, where 
$\tau(n)$ is the $n$-th Fourier coefficient of the unique normalized cusp form of weight $12$ and full level, and 
the prime $691$ appears in the numerator of $\zeta(12)/\pi^{12}$ for the Riemann zeta function $\zeta(s)$. 
Searching for such congruences, it is shown that 
the prime $67$ appears in the numerator of $L(6,\chi)/(\pi^6 \sqrt{5})$, where $\chi$ is the unique nontrivial quadratic Dirichlet character modulo $5$ and 
$L(s,\chi)$ its Dirichlet $L$-function, giving rise to a congruence 
$f_\chi \ \equiv \ E^\circ_{6, \chi} \pmod{67}$
between a cusp form $f_\chi$ and an Eisenstein series $E^\circ_{6, \chi}$ of weight $6$ on $\Gamma_0(5)$ with nebentypus character $\chi.$ 

\end{abstract}

\maketitle

 \def\R{\mathbb{R}}
\def\C{\mathbb{C}}
\def\Z{\mathbb{Z}}
\def\Q{\mathbb{Q}}
\def\A{\mathbb{A}}
\def\F{\mathbb{F}}
 \def\J{\mathbb{J}}
\newcommand\D{{\mathbb{ D}}}
 \def\L{\mathbb{L}}
\def\bG{\mathbb{G}}
\def\N{\mathbb{N}}
\def\BH{\mathbb{H}}
\newcommand\Qi{\mathbb{Q}({\bf i})}

\def\bfG{\mathbf{G}}
\def\bfB{\mathbf{B}}
\def\bfT{\mathbf{T}}
\def\bfU{\mathbf{U}}
\def\bfP{\mathbf{P}}
\def\bfM{\mathbf{M}}
\def\bfN{\mathbf{N}}
\def\bfA{\mathbf{A}}

\def\LN{{}^L\!N}
\def\LP{{}^L\!P}
\def\LM{{}^L\!M}
\def\LA{{}^L\!A}
\def\LG{{}^L\!G}
\def\Ln{{}^L\mathfrak{n}}
\def\Lg{{}^L\mathfrak{g}}

\newcommand{\vless}{\rotatebox[origin=c]{-90}{$<$}}
\newcommand{\vgreat}{\rotatebox[origin=c]{90}{$<$}}
\newcommand{\vgreater}{\rotatebox[origin=c]{90}{$\leq$}}

\newcommand\Dm{\D_\lambda}
\newcommand\Dmp{\D_{\lambda^\prime}}
 \newcommand\Dim{\D_{\io\lambda}}
 \newcommand\Dimp{\D_{\io\lambda^\prime}}
\newcommand\Dum{\D_{\ul{\lambda}}}
\newcommand\Dium{\D_{\io\ul{\lambda}}} 
\newcommand\DumN{\D_{\ul{\lambda}-N\gamma_P}} 
\newcommand\DiumN{\D_{\io\ul{\lambda}-N\gamma_P}} 
\newcommand\cal{\mathcal}
\newcommand\SMK{{\cal S}^M_{K^M_f}}
\newcommand\tMZl{\tM_{\lambda,\Z}} 
\newcommand\Gm{{\mathbb G}_m}
\newcommand\cA{\cal A}
\newcommand\cC{\cal C}
\newcommand\calL{\cal L}
\newcommand\cO{\cal O}
\newcommand\cU{\cal U}
\newcommand\cK{\cal K}   
\newcommand\cW{\cal W}     
\newcommand\HH{{\cal H}}
\newcommand\cF{\mathcal{F}} 
\newcommand\G{\mathcal{G}}
\newcommand\cB{\mathcal{B}}
\newcommand\cT{\mathcal{T}}
\newcommand\cS{\mathcal{S}}
\newcommand\cP{\mathcal{P}}
\newcommand\Exp{\mathcal{E}xp}

\newcommand\GL{{ \rm  GL}}
\newcommand\Gl{{ \rm  GL}}
\newcommand\U{{ \rm  U}}
\def\SU{{\rm SU}}
\def\S{{\bf S}}
\newcommand\Gsp{{\rm Gsp}}
\newcommand\Lie {{ \rm Lie}} 
\newcommand\Sl{{ \rm  SL}}
\newcommand\SL{{ \rm  SL}}
\newcommand\SO{{ \rm  SO}}
\newcommand\rO{{\rm  O}}
\newcommand\GU{{\rm  GU}}
\newcommand{\Sp}{{\rm Sp}}
\newcommand\Ad{{\rm Ad}}
\newcommand\Sym{{\rm Sym}}

\def\ringO{\mathcal{O}}
\def\idealP{\mathfrak{P}} 
\def\g{\mathfrak{g}}
\def\a{\mathfrak{a}}
\def\k{\mathfrak{k}}
\def\z{\mathfrak{z}}
\def\m{\mathfrak{m}}
\def\s{\mathfrak{s}}
\def\c{\mathfrak{c}}
\def\b{\mathfrak{b}}
\def\t{\mathfrak{t}}
\def\q{\mathfrak{q}}
\def\l{\mathfrak{l}}
\def\gl{\mathfrak{gl}}
\def\sl{\mathfrak{sl}}
\def\u{\mathfrak{u}} 
\def\fp{\mathfrak{p}} 
\def\p{\mathfrak{p}}   
\def\r{\mathfrak{r}}
\def\fd{\mathfrak{d}}
\def\fR{\mathfrak{R}}
\def\fI{\mathfrak{I}}
\def\fJ{\mathfrak{J}}
\def\i{\mathfrak{i}}
\def\perm{\mathfrak{S}}
\newcommand\fg{\mathfrak g}
\newcommand\fk{\mathfrak k}
\newcommand\fgK{(\mathfrak{g},K_\infty^0)}
\newcommand\gK{ \mathfrak{g},K_\infty^0 }

\newcommand\ul{\underline} 
\newcommand\tp{{  {\pi}_f}}
\newcommand\tv{ {\pi}_v}
\newcommand\ts{{ {\pi}_f}}
\newcommand\pts{ {\pi}^\prime_f}
\newcommand\usf{\ul{\pi}_f}
\newcommand\pusf{\ul{\pi}^\prime_f}
\newcommand\usvp{\ul{\pi}^\prime_v}
\newcommand\usv{\ul{\pi}_v}
\newcommand\io{{}^\iota}
\newcommand\uls{{\underline{\pi}}}

\newcommand\Spec{\hbox{\rm Spec}} 
\newcommand\SGK{\mathcal{S}^G_{K_f}}
\newcommand\SMP{\mathcal{S}^{M_P}}
 \newcommand\SGn{\mathcal{S}^{G_n}}
 \newcommand\SGp{\mathcal{S}^{G_{n^\prime}}}
\newcommand\SMPK{\mathcal{S}^{M_P}_{K_f^{M_P}}}
\newcommand\SMQ{\mathcal{S}^{M_Q}}
\newcommand\SMp{\mathcal{S}^{M }_{K_f^M}}
\newcommand\SMq{\mathcal{S}^{M^\prime}_{K_f^{M^\prime}}}
\newcommand\uSMP{\ul{\mathcal{S}}^{M_P}}
\newcommand\SGnK{\mathcal{S}^{G_n}_{K_f}}
\newcommand\SG{\mathcal{S}^G}
\newcommand\SGKp{\mathcal{S}^G_{K^\prime_f}}
\newcommand\piKK{{ \pi_{K_f^\prime,K_f}}}
\newcommand\piKKpkt{\pi^{\pkt}_{K_f^\prime,K_f}}
\newcommand\BSC{ \bar{\mathcal{S}}^G_{K_f}}
\newcommand\PBSC{\partial\SGK}
\newcommand\pBSC{\partial\SG}
\newcommand\PPBSC{\partial_P\SGK}
\newcommand\PQBSC{\partial_Q\SGK}
\newcommand\ppBSC{\partial_P\mathcal{S}^G}
\newcommand\pqBSC{\partial_Q\mathcal{S}^G}
\newcommand\prBSC{\partial_R\mathcal{S}^G}
\newcommand \bs{\backslash} 
 \newcommand \tr{\hbox{\rm tr}}
 \newcommand\ord{\text{ord}}
\newcommand \Tr{\hbox{\rm Tr}}
\newcommand\HK{\mathcal{H}^G_{K_f}}
\newcommand\HKS{\mathcal{H}^G_{K_f,\place}}
\newcommand\HKv{\mathcal{H}^G_{K_v}}
\newcommand\HGS{\mathcal{H}^{G,\place}}
\newcommand\HKp{\mathcal{H}^G_{K_p}}
\newcommand\HKpo{\mathcal{H}^G_{K_p^0}}
\newcommand\ch{{\bf ch}}

\newcommand\M{\mathcal{M}}
\newcommand\Ml{\M_\lambda}
\newcommand\tMl{\tilde{\Ml}}
\newcommand\tM{\widetilde{\mathcal{M}}}
\newcommand\tMZ{\tM_\Z}
\newcommand\tsigma{\ul{\pi}}
\newcommand \pkt{\bullet}
\newcommand\tH{\widetilde{\mathcal{H}}}
\newcommand\Mot{{\bf M}} 
\newcommand\eff{{\rm eff}}
\newcommand\Aql{A_{\q}(\lambda)}
\newcommand\wl{w\cdot\lambda}
\newcommand\wlp{w^\prime\cdot\lambda} 

\def\w{{\bf w}} 
\def\d{{\sf d}}
\def\e{{\bf e}} 
\def\x{{\tt x}}
\def\y{{\tt y}}
\def\v{{\sf v}}
\def\q{{\sf q}} 
\def\ff{{\bf f}}
\def\bk{{\bf k}}
 
\def\Ext{{\rm Ext}}
\def\Aut{{\rm Aut}}
\def\Hom{{\rm Hom}}
\def\Ind{{\rm Ind}}
\def\Asai{{\rm Asai}}
\def\aInd{{}^{\rm a}{\rm Ind}}
\def\aIndPG{\aInd_{\pi_0(P(\R)) \times P(\A_f)}^{\pi_0(G(\R)) \times G(\A_f)}}
\def\aIndQG{\aInd_{\pi_0(Q(\R)) \times Q(\A_f)}^{\pi_0(G(\R)) \times G(\A_f)}}
\def\Gal{{\rm Gal}}
\def\End{{\rm End}} 
\newcommand\Coh{{\rm Coh}}  
\newcommand\Eis{{\rm Eis}}
\newcommand\Res{\mathrm{Res}}
\newcommand\place{\mathsf{S}}
\newcommand\emb{\mathcal{I}} 
\newcommand\LB{\mathcal{L}}  
\def\Hod{{\mathcal{H}od}}
\def\Crit{{\rm Crit}}

\def\bfpi{\mathbf{\Pi}}
\def\bfdelta{\mathbf{\Delta}}

\section{Introduction}
\label{sec:intro}
In a landmark paper published in 1916 Srinivasa Ramanujan \cite{ramanujan} studied the function: 
$$
\Delta(z) \ := \ q \prod_{n=1}^\infty (1-q^n)^{24} \ := \ \sum_{n=1}^\infty \tau(n) q^n, \quad q := e^{2\pi i z}. 
$$
Ramanujan discovered the remarkable congruence that for any prime $p$ one has:
$$
\tau(p) \ \equiv \ p^{11} + 1 \pmod{691}.
$$
The prime $691$, which may be called Ramanujan's congruence prime, also appears elsewhere in number theory as being an irregular prime in the sense of Kummer because 
it divides the numerator of the Bernoulli number $B_{12}.$ For this latter reason, the prime $691$ appears in the numerator of the rational number $\zeta(12)/\pi^{12}$, where $\zeta(s)$ is the Riemann $\zeta$-function; the rationality of $\zeta(12)/\pi^{12}$ is a classical result due to Euler from 1730's. 
Ramanujan's congruence can also be interpreted as the congruence 
$$
\Delta \ \equiv \ E_{12} \pmod{691}, 
$$
where $\Delta$ is the unique (up to scaling) cusp form of weight $12$ on the modular group $\SL_2(\Z),$ and $E_{12}$ is (the normalized) Eisenstein series of weight $12$ on
$\SL_2(\Z).$ See Manin \cite{manin} for the above congruences, and for many more such congruences, for example, for weight $16$ modular forms modulo the prime $3617$ because $3617$ appears in the numerator of $\zeta(16)/\pi^{16}$. 


The purpose of this article is to indicate a computational search for other Ramanujan's congruence primes (such as $691$ and $3617$) by looking at the numerators of the special values of $L$-functions attached to Dirichlet characters and provide evidence for the existence of such congruences. For example, the prime $67$ appears in the numerator of $L(6,\chi)/(\pi^6 \sqrt{5})$, where $\chi$ is the unique nontrivial quadratic Dirichlet character modulo $5$, suggesting the existence of a Ramanujan like congruence (see Remark~\ref{rem:r-prime}). 
Indeed, we prove that there is a congruence 
$$
f_\chi \ \equiv \ E^\circ_{6, \chi} \pmod{67}, 
$$
where $f_\chi$ is a cusp form and $E^\circ_{6, \chi}$ an Eisenstein series (defined in Section~\ref{sec:hecke-eisenstein}) both of weight $6$, on the Hecke congruence subgroup $\Gamma_0(5)$, and with nebentypus character $\chi.$

\bigskip
\section{Special values of the Riemann $\zeta$-function}

A famous result of Euler states that the values of the Riemann zeta function at even positive integers are given by: 
$$
\zeta(2m) \ = \ (-1)^{m-1} \frac{(2\pi)^{2m}}{2 (2m)!} B_{2m},
$$
where the Bernoulli numbers $B_m$ are defined by the power series: 
$$
\frac{t e^t}{e^t - 1} \ = \ \sum_{n=0}^\infty B_m \frac{t^m}{m!}. \quad (\mbox{Note: $B_m \in \Q$}.)
$$
(See, for example, Neukirch \cite[Sec.\,VII.1]{neukirch}.) For brevity, let $Z_{2m} = \zeta(2m)/\pi^{2m},$ and write the rational number $Z_{2m}$ as $N_{2m}/D_{2m}$, for relatively prime integers $N_{2m}$ and $D_{2m}$; these values for $2m \in \{2,4,6,\dots, 20\}$ are: 
$$
\begin{array}{|c|c|c|}
\hline
2m & N_{2m} & D_{2m}  \\
\hline \hline
2 & 1 & 2 \cdot 3 \\
\hline 
4 & 1 & 2 \cdot 3^2 \cdot 5 \\
\hline 
6 & 1 & 3^3 \cdot 5 \cdot 7 \\
\hline 
8 & 1 & 2 \cdot 3^3 \cdot 5^2 \cdot 7 \\
\hline 
10 & 1 & 3^5 \cdot 5 \cdot 7 \cdot 11 \\
\hline 
12 & 691 & 3^6 \cdot 5^3 \cdot 7^2 \cdot 11 \cdot 13 \\
\hline 
14 & 2 & 3^6 \cdot 5^2 \cdot 7 \cdot 11 \cdot 13 \\
\hline 
16 & 3617 & 2 \cdot 3^7 \cdot 5^4 \cdot 7^2 \cdot 11 \cdot 13 \cdot 17 \\
\hline 
18 & 43867 & 3^9 \cdot 5^3 \cdot 7^3 \cdot 11 \cdot 13 \cdot 17 \cdot 19 \\
\hline 
20 & 283 \cdot 617 & 3^9 \cdot 5^5 \cdot 7^2 \cdot 11^2 \cdot 13 \cdot 17 \cdot 19 \\
\hline
\end{array}
$$

\smallskip
The denominators $D_{2m}$ have small primes appearing in them. Typically, most of the primes in the range $1$ through $2m$ show up because of the term $(2m)!$ in Euler's theorem; 
also if we use the Clausen - von Staudt theorem (see, for example, \cite[Thm.\,3.1]{arakawa-ibukiyama-kaneko}) which states that 
$$
B_n + \sum_{(p-1)|n} \frac{1}{p} \in \Z, 
$$
then whenever $2m+1$ is a prime, 
it appears in $D_{2m}$. For example, if $2m \in \{2, 4, 6, 10, 12, 16, 18\}$ then $2m+1$ is a prime that appears in $D_{2m}$ in the table above.  

\smallskip

The numerator is a different story! Occasionally, some large prime appears in the numerator. The first such instance is $691$ appearing in the numerator of $Z_{12}.$ 
Such strange primes are called Ramanujan's congruence primes. This is not a rigorous definition, because one is not defining the meaning of `strange', however, it constitutes a useful working principle, which is the case being made in the article. One may attempt 
a rigorous definition by saying that a prime $p > 2m+1$ appearing in the numerator of $Z_{2m}$ may be called a Ramanujan's congruence prime. For example, the prime $3617$ appears in the numerator of $\zeta(16)$ is indeed a Ramanujan congruence prime. Suppose 
$f(z) = \sum_{n=1}^\infty a_n q^n $
is the unique cusp form of weight $16$ for $\SL_2(\Z)$ normalized as $a_1 = 1$, 
then it is an easy consequence of Manin's `coefficients theorem' \cite[Thm.\,I.3]{manin} that 
$$
a_p \equiv p^{15} + 1 \pmod{3617}.
$$
See the table on p.\,383 of \cite{manin} for this congruence and other such examples. Whether every prime larger than $2m+1$ appearing in the numerator of $Z_{2m}$ 
is indeed a congruence prime is not clear to us. 

\bigskip
\section{Special values of Dirichlet $L$-functions}

Fix an odd prime $p$, and let $\chi_p$ denote the unique nontrivial quadratic Dirichlet character modulo $p$ given by the Legendre symbol $\left( \frac{\cdot}{ p}\right)$. 
For notational brevity, we will often denote $\chi_p$ simply by $\chi$, suggestive of the fact that most of what follows makes sense for a general Dirichlet character $\chi$; 
however, we will exclusively be working with our quadratic $\chi = \chi_p$. 
The $L$-function attached to $\chi$ is the Dirichlet series: 
$$
L(s, \chi) \ = \ \sum_{n=1}^\infty \frac{\chi(n)}{n^s}. 
$$
Of course, one may enlarge the context to the usual framework of Dirichlet characters modulo any positive integer; but already in the simplest situation of 
a quadratic character modulo an odd prime $p,$ one sees interesting new congruences. In this short note we present one such illustrative example. 
If $p=5$, then the $L$-function has the following shape: 
$$
L(s, \chi) \ = \  \frac{1}{1^s} - \frac{1}{2^s} - \frac{1}{3^s} + \frac{1}{4^s} + \frac{1}{6^s}  - \frac{1}{7^s} - \frac{1}{8^s} + \frac{1}{9^s} + \frac{1}{11^s} \cdots  
$$

Henceforth, assume also that $p \equiv 1 \pmod{4}.$ This has the consequence that $\chi(-1) = 1.$ The following result is a generalization of Euler's formula
(see Neukirch \cite[Cor.\,VII.2.10]{neukirch}). 

\begin{thm}[Leopoldt (1958)]
Let $\chi$ denote the unique nontrivial quadratic Dirichlet character modulo an odd prime $p\equiv 1 \pmod{4}.$ 
For any positive integer $m$ we have:
$$
L(2m, \chi) \ = \ (-1)^{m-1}\, \frac{\sqrt{p}}{2} \left(\frac{2\pi}{p}\right)^{2m} 
\frac{B_{2m, \chi}}{(2m)!}. 
$$
\end{thm}

The term $\sqrt{p}$ is the Gauss sum of $\chi$, and the numbers $B_{2m,\chi}$ are the generalized Bernoulli numbers introduced by Leopoldt by 
the generating function:
$$
\sum_{a = 1}^{p} \chi(a) \frac{t e^{at}}{e^{pt} - 1} \ = \ 
\sum_{n=0}^\infty B_{n,\chi} \frac{t^n}{n!}.
$$

Under the above hypotheses on $\chi$, it follows that 
$B_{0,\chi} = 0$
and $B_{2k+1,\chi} = 0$ for all $k \geq 0.$ 
One has the well-known integrality result of Carlitz \cite{carlitz} that bounds the denominator: 

\begin{thm}[Carlitz]
Let $\chi$ denote the unique nontrivial quadratic Dirichlet character modulo an odd prime $p\equiv 1 \pmod{4}.$ 
Then $p \cdot B_{2m,\chi}$ is an integer for all integers $m \geq 1.$
\end{thm}

\begin{rem}
\label{rem:r-prime}{\rm 
It follows from these theorems of Leopoldt and Carlitz that the numerator of the rational number
$$
\frac{L(2m, \chi)}{\pi^{2m} \sqrt{p}}
$$ 
is essentially the same as the integer $p \cdot B_{2m,\chi}$. In particular, 
to search for Ramanujan type congruence primes, one can search for `strange' or suitably large primes appearing in the integer $p \cdot B_{2m,\chi}$. As a working principle, 
we may take any prime $\ell > \max\{p, 2m+1\}$ that divides the integer $p \cdot B_{2m,\chi}$ as a candidate for Ramanujan's congruence prime for $\chi$ and $2m.$
}\end{rem}

\subsection{Computing $B_{2m, \chi}$}

Computing generalized Bernoulli numbers is discussed in detail in Stein's book \cite[Sec.\,5.2]{stein}. Since we are are working in a simplistic situation, we can also 
compute these numbers using some basic information which we now review.

\subsubsection{Computing $B_{2m, \chi}$ via Bernoulli polynomials}

The Bernoulli polynomial $B_n(x)$ is a monic polynomial of degree $n$ in the variable $x$ with coefficients in $\Q$ given by the generating series: 
$$
\frac{te^{xt}}{e^t-1}  \ = \ 
 \sum_{n=0}^\infty B_n(x) \frac{t^n}{n!}. 
$$
In terms of the Bernoulli numbers we have (\cite[Prop.\,4.9]{arakawa-ibukiyama-kaneko}):  
$$
B_n(x) \ = \ \sum_{j=0}^n (-1)^j \binom{n}{j} B_j x^{n-j}.
$$
The generalized Bernoulli numbers are given via the formula (\cite[(4.1)]{arakawa-ibukiyama-kaneko}): 
\begin{equation}
\label{eqn:gen_bernoulli_poly}
B_{n,\chi} \ = \ p^{n-1} \sum_{a=1}^p \chi(a) B_n(a/p).
\end{equation}
In other words, if we know the Bernoulli numbers $B_n$ then we can compute the generalized Bernoulli numbers $B_{n,\chi}.$

\subsubsection{Computing $B_{2m, \chi}$ via recursion}

For $\chi$, a Dirichlet character modulo $p$, we define
$$
S_{\chi}(n) \ = \ \sum_{a=1}^{p-1} \chi(a) a^n. 
$$

The following recursive formula, possibly due to Carlitz, is well-known (see, for example, Agoh \cite{agoh}): 

\begin{thm}
\label{thm:gen_bernoulli_recursion}
Let $\chi$ denote the unique nontrivial quadratic Dirichlet character modulo an odd prime $p\equiv 1 \pmod{4}.$ 
For each even integer $2m$, we have 
$p \cdot B_{2, \chi} = S_\chi(2),$ and for $m \geq 2$ we have: 
$$
B_{2m, \chi} = \frac{1}{p} \left[S_\chi(2m) - \left(\sum_{j=0}^{m-2} \binom{2m}{2j+1} 
\frac{B_{2j+2,\chi}}{2j+2} 
p^{2m-2j-1} \right)\right].
$$
\end{thm}

\subsection{Values of $B_{2m, \chi}$ for $p=5$}

We have the following table of values $p \cdot B_{2m,\chi}$ for $p=5$: 
$$
\begin{array}{|c|c|c|}
\hline
2m & p \cdot B_{2m,\chi}  & \mbox{prime factorization} \\
\hline \hline
2 & 4  & 2^2\\
\hline 
4 & -40 & -1 \cdot 2^3 \cdot 5 \\
\hline 
6 & 804  & 2^2 \cdot 3 \cdot 67\\
\hline 
8 & -28880 & -1 \cdot 2^4 \cdot 5 \cdot 19^2\\
\hline 
10 & 1651004 & 2^2 \cdot 191 \cdot 2161\\
\hline 
12 & -138110520 & -1 \cdot 2^3 \cdot 3 \cdot 5 \cdot 1150921\\
\hline 
14 & 15920571604 & 2^2 \cdot 7 \cdot 17 \cdot 33446579 \\
\hline 
16 & -2419747948960 & -1 \cdot  2^5 \cdot  5 \cdot 457 \cdot  33092833 \\
\hline 
18 & 468896302250604 & 2^2 \cdot  3^2 \cdot 41 \cdot  317680421579 \\
\hline 
20 & -112834502909928192 & -1 \cdot 2^8 \cdot 3 \cdot 146919925663969 \\
\hline
\end{array}
$$

\medskip
\begin{rem}{\rm
A word of caution in computing these numbers. Intitially, we wrote a Python code using Theorem~\ref{thm:gen_bernoulli_recursion}. Interestingly, the code gave a non-integral value 
for $5 \cdot B_{16,\chi}$. Then we wrote a Python code using \eqref{eqn:gen_bernoulli_poly} and it was giving a nonzero value for $B_{17,\chi}$. The answer is that Python does not handle large numbers very well. The above table of values is computed on SAGE, and one gets the same list of numbers computed by either of the two methods. 
}\end{rem}

\bigskip
\section{Ramanujan type congruence in $M_6(\Gamma_0(5), \chi)$ modulo $67$} 

If $\chi$ is the quadratic Dirichlet character modulo $p=5$, then the prime $67$ divides $5 \cdot B_{6,\chi}$. In this section we verify that $67$ is indeed a Ramanujan's congruence prime. Towards this one looks at the space $M_6(\Gamma_0(5), \chi)$ of modular forms for weight $6$, level $5$, i.e., for the Hecke congruence subgroup $\Gamma_0(5)$, and with nebentypus $\chi$.

\bigskip
\subsection{The Hecke--Eisenstein series $E^{\circ }_{6, \chi} \in M_6(\Gamma_0(5), \chi)$}
\label{sec:hecke-eisenstein}
For the moment, let $\chi$ be the unique nontrivial quadratic Dirichlet character modulo an odd prime $p\equiv 1 \pmod{4}.$ Let $k = 2m$ be an even positive integer. Hecke 
\cite{hecke} studied a family of Eisenstein series; we need one such Eisenstein series--our notation is partly adapted from Shimura \cite[Sec.\,2]{shimura}). Let 
$$
E_{k, \chi}(z) \ := \ \sum_{(m,n) \in \Z^2 \setminus \{(0,0)\}} \frac{\chi(n)}{(mpz + n)^k}, 
$$
as a function of a complex variable $z$ in the upper halfplane. 
Hecke proved that $E_{k, \chi}$ is an element of $M_k(\Gamma_0(p), \chi)$, and has the Fourier expansion: 
$$
E_{k, \chi}(z) \ = \ 2 L(k, \chi) + \frac{2 \sqrt{p} (-2\pi i)^k}{p^k (k-1)!} \sum_{n=1}^\infty \left\{\sum_{d|n} \chi(d) d^{k-1} \right\} q^n, \quad (q = e^{2\pi i z}).
$$
(See, for example, Shimura \cite[(3.4)]{shimura}.) Define a normalized Eisenstein series by
$$
E^\circ_{k, \chi}(z) \ := \ \left(\frac{2 \sqrt{p} (-2\pi i)^k}{p^k (k-1)!}\right)^{-1} E_{k, \chi}(z),
$$
and apply Leopoldt's theorem to the constant term, to get: 
$$
E^\circ_{k, \chi}(z) \ = \ - \frac{B_{k,\chi}}{2k} +  \sum_{n=1}^\infty \left\{\sum_{d|n} \chi(d) d^{k-1} \right\} q^n.
$$
(See, for example, Stein \cite[(5.3.1)]{stein}.) In the special case of $p=5$ and $k = 6$, 
we get the following $q$-expansion of the normalized Eisenstein series: 
\begin{multline}
\label{eqn:e-chi}
    E^\circ_{6, \chi}(z) = -\frac{67}{5} +  q - 31 q^{2} - 242 q^{3} + 993 q^{4} +  q^{5} + 7502 q^{6} \\ 
    - 16806 q^{7} - 31775 q^{8} + 58807 q^{9} - 31 q^{10} + O(q^{11}). 
\end{multline}

\bigskip
\subsection{The cusp form $f_\chi \in S_6(\Gamma_0(5), \chi)$}
The space $S_6(\Gamma_0(5), \chi)$ of cusp forms in $M_6(\Gamma_0(5), \chi)$ is two-dimensional having as basis eigenforms with coefficients in $\Q(\alpha)$, 
where $\alpha = -\sqrt{-44}$ is a root of the polynomial $x^2 + 44$. The two eigenforms are conjugate by $\Gal(\Q(\alpha)/\Q)$. As usual, we take them normalized by the requirement that the 
first Fourier coefficient is $1$. One of these forms, denoted (say) $f_\chi$, has the $q$-expansion:
\begin{multline}
\label{eqn:f_chi}
 f_\chi(z)  = q + \alpha q^{2} - 3 \alpha q^{3} - 12 q^{4} + \left(5 \alpha - 45\right) q^{5} + 132 q^{6} \\ 
 - 9 \alpha q^{7}  + 20 \alpha q^{8} - 153 q^{9} + \left(-45 \alpha - 220\right) q^{10}  + O(q^{11}) . 
\end{multline}
This $f_\chi$ is the nontrivial $\Gal(\Q(\alpha)/\Q)$-conjugate of the form denoted 5.6.b.a in LMFDB \cite{lmfdb}.

\bigskip
\subsection{Congruence mod $67$}
The prime $67$ splits in the extension $\Q(\alpha)$, since the Legendre symbol $(-44/67) = 1.$ Suppose $67$ factors as $\p_1 \p_2$; to check if $\p_j$ divides 
an element $x \in \Q(\alpha)$, it is simpler to check if $67$ divides the norm $N_{\Q(\alpha)/\Q}(x)$ of $x$. The factorization of the norm of the difference $a_n(f_\chi) - a_n(E^\circ_{6, \chi})$
of the $n$-th Fourier coefficients: 
\begin{equation}\label{table:congruence}
 \begin{tabular}{|l|l|} \hline
$n$ & $N_{\mathbb{Q}(\alpha)/\mathbb{Q}}\left(a_n(f_\chi) - a_n(E^\circ_{6, \chi}) \right)$ \\ \hline \hline
$0$ & $5^{-1} \cdot 67$ \\ \hline
$1$ & $0$ \\ \hline
$2$ & $3 \cdot 5 \cdot 67$ \\ \hline
$3$ & $2^{4} \cdot 5 \cdot 11 \cdot 67$ \\ \hline
$4$ & $3^{2} \cdot 5^{2} \cdot 67^{2}$ \\ \hline
$5$ & $2^{4} \cdot 3 \cdot 67$ \\ \hline
$6$ & $2^{2} \cdot 5^{2} \cdot 11^{2} \cdot 67^{2}$ \\ \hline
$7$ & $2^{4} \cdot 3^{2} \cdot 5^{2} \cdot 67 \cdot 1171$ \\ \hline
$8$ & $3 \cdot 5^{2} \cdot 67 \cdot 200929$ \\ \hline
$9$ & $2^{8} \cdot 5^{2} \cdot 11^{2} \cdot 67^{2}$ \\ \hline
$10$ & $3^{4} \cdot 23 \cdot 67$ \\ \hline
\end{tabular}
\end{equation}

\medskip
The following criterion due to Sturm \cite{sturm} is extremely useful in verifying congruences between two modular forms. (For more details see Stein's book \cite[Sec.\,9.4]{stein}.)

\begin{thm}[Sturm]
Let $f,g \in M_k(\Gamma_0(N), \omega, \cO)$ be holomorphic modular forms of weight $k$, level $N$, nebentypus $\omega$, with all the Fourier coefficients $a_n(f)$ and $a_n(g)$ 
in the ring of integers 
$\cO$ of a number field. Suppose $\m$ is an integral ideal of $\cO$ such that 
$$
a_n(f) \ \equiv \ a_n(g) \pmod{\m}, \quad \mbox{for all $n < \frac{k \cdot [\SL_2(\Z): \Gamma_0(N)]}{12}$}, 
$$
then $f \equiv g \pmod{\m},$ i.e., $a_n(f)  \equiv  a_n(g) \pmod{\m}$ holds for all $n$. 
\end{thm}

Apply this to $M_6(\Gamma_0(5), \chi);$  it suffices to check that the congruence holds for the $n$-th Fourier coefficient for 
$n < \frac{6 \cdot [\SL_2(\Z): \Gamma_0(5)]}{12} = 3$. From the table \eqref{table:congruence}, we deduce the following result: 

\medskip
\begin{thm}
\label{thm:main-example}
Let $f_\chi \in S_6(\Gamma_0(5), \chi)$ be the cusp form with Fourier expansion \eqref{eqn:f_chi}, 
and $E^\circ_{6, \chi} \in M_6(\Gamma_0(5), \chi)$ the normalized Hecke--Eisenstein series with Fourier expansion
\eqref{eqn:e-chi}. Then we have the congruence: 
 $$
f_\chi \ \equiv \ E^\circ_{6, \chi} \pmod{67}.
 $$
\end{thm}

\section{Some comments and acknowledgements} 

First and foremost we would like to acknowledge the influence of G\"unter Harder's ideas on congruences and the special values of $L$-functions. 
In particular, the second author learnt the idea of looking for such strange primes from conversations with Harder, especially concerning the ideas around his conjecture on congruences between elliptic and Siegel modular forms \cite{harder-123}. It is well-known to experts in the cohomology of arithmetic groups, that the theory of Eisenstein cohomology pioneered by Harder provides a theoretical framework to prove congruences as in Theorem~\ref{thm:main-example}. In his book \cite{harder-book}, Harder addresses the problem of understanding the denominators of Eisenstein classes; for example, Ramanujan's original congruence prime $691$ appears in the denominator of an Eisenstein class in the cohomology of the upper half-plane modulo $\SL_2(\Z)$ with coefficients in the sheaf of $\Q$-vector spaces given by the unique irreducible representation $\Sym^{10}(\Q^2)$ of $\SL_2/\Q$ of dimension $11$. 

\medskip

Let $X_0(5) = \Gamma_0(5)\backslash \HH$ be the upper half-plane $\HH$ modulo $\Gamma_0(5)$. Let $V_{6,\chi}$ denote the $5$-dimensional irreducible representation 
$\Sym^{4}(\Q^2) \otimes \chi$ of $\Gamma_0(5)$, and $\widetilde{V}_{6,\chi}$ the corresponding sheaf of $\Q$-vector spaces on $X_0(5)$. We expect to find the prime $67$ to appear in the 
denominator of an Eisenstein cohomology class in $H^1(X_0(5), \widetilde{V}_{6,\chi})$. Note that $V_{6,\chi}$, because of the twisting by $\chi$, is not  an algebraic representation 
of the algebraic group $\SL(2)/\Q$; nevertheless, such a representation and the corresponding sheaf are considered in the context of the Eichler--Shimura isomorphism; 
see Shimura \cite[Sec.\,8.2]{shimura}.
More generally, 
if $\chi$ is the nontrivial quadratic Dirichlet character modulo an odd prime $p \equiv 1 \pmod{4}$, and $\ell$ is a `Ramanujan congruence prime' appearing in the numerator of 
$L(2m, \chi)/(\pi^{2m}\sqrt{p})$, and now letting $V_{2m,\chi} = \Sym^{2m-2}(\Q^2) \otimes \chi$, we would expect $\ell$ to appear as a denominator
of an Eisenstein cohomology class in $H^1(X_0(p), \widetilde{V}_{2m,\chi})$. In our opinion, this example, and its obvious generalizations, should be studied further.   

\medskip

Ramanujan's congruence and its generalizations is a well-known theme. For example, there is a generalization due to Gaba and Popa \cite{gaba-popa}, that subsumes many previous generalizations. Let us add that the congruence in Theorem \ref{thm:main-example} is different from the main results of \cite{gaba-popa}. The reader is also referred to 
Swinnerton-Dyer \cite{s-dyer} for the deep connections between these congruences in the theory of modular forms with Galois representations.

\bigskip
\bigskip

\noindent
{\Small \textit{Acknowledgements:} 
We thank Cris Poor for a helpful conversation on computing generalized Bernoulli numbers and drawing our attention 
to the beautiful book \cite{arakawa-ibukiyama-kaneko}. We are deeply grateful to P. Narayanan, currently a doctoral student at IISER Pune, India, 
who provided us with the details verifying the congruence modulo $67$. We thank Siddhartha Mitra, currently a graduate student in computational neuroscience at SUNY, Brooklyn, who clarified some interesting computational features with Python and SAGE. Thanks also to Ken Ono for pointing us to Sturm \cite{sturm} with the useful criterion to verify congruences. Finally, 
we thank the referee for several suggestions to clean up our exposition.}

\bigskip
\bigskip

\bibliographystyle{plain}

\end{document}